\documentclass[10pt]{article}
\usepackage{amsmath,amssymb,amsthm, hyperref, euscript}
\usepackage[matrix,arrow,curve]{xy}
\usepackage{graphicx}
\usepackage{tabularx}
\usepackage{float}
\usepackage{hyperref}

\usepackage[T1]{fontenc}

\usepackage[sc]{mathpazo}
\usepackage[usenames]{color}

\DeclareFontFamily{T1}{pzc}{}
\DeclareFontShape{T1}{pzc}{m}{it}{1.8 <-> pzcmi8t}{}
\DeclareMathAlphabet{\mathpzc}{T1}{pzc}{m}{it}

\title{Noncommutative Generalization of Wilson Lines}

\theoremstyle{plain}
\newtheorem{prop}{Proposition}[section]

\newtheorem{lem}[prop]{Lemma}

\theoremstyle{definition}
\newtheorem{defn}[prop]{Definition}
\newtheorem{empt}[prop]{}
\newtheorem{exm}[prop]{Example}
\newtheorem{rem}[prop]{Remark}

\theoremstyle{remark}

\chardef\bslash=`\\ 

\newcommand{\Om}{\Omega}

\newcommand{\lb}{\lambda}

\newcommand{\rar}{\rightarrow}

\newbox\ncintdbox \newbox\ncinttbox 
\setbox0=\hbox{$-$} \setbox2=\hbox{$\displaystyle\int$}
\setbox\ncintdbox=\hbox{\rlap{\hbox
    to \wd2{\hskip-.125em \box2\relax\hfil}}\box0\kern.1em}
\setbox0=\hbox{$\vcenter{\hrule width 4pt}$}
\setbox2=\hbox{$\textstyle\int$} \setbox\ncinttbox=\hbox{\rlap{\hbox
    to \wd2{\hskip-.175em \box2\relax\hfil}}\box0\kern.1em}

\begin{document}
\maketitle  \setlength{\parindent}{0pt}
\begin{center}
\author{}
{\textbf{Petr R. Ivankov*}\\
e-mail: * monster.ivankov@gmail.com \\
}
\end{center}

\vspace{1 in}

\begin{abstract}
\noindent

A classical Wilson line is a correspondence between closed paths and elements of a gauge group. However the noncommutative geometry does not have closed paths. But noncommutative geometry have good generalizations of both: the covering projection, and the group of covering transformations. These notions are used for a construction of noncommutative Wilson lines. Wilson lines can also be constructed as global pure gauge fields on the universal covering space. The noncommutative analog of this construction is also developed.

\end{abstract}
\tableofcontents

\section*{Foreword}
\paragraph{}At 1970s I had a strong wish for understanding of physics. I had been reading many books written by physicists, but I did not understand them. Later I started to study math and then I became to understand physics. General relativity became very clear after knowledge of Riemann geometry. Physicist uses the "vector" notion, but  mathematician uses rigorous  notion of a vector space. It is clear what space mathematicians mean: Banach, Hilbert or Fr\'{e}chet. A gauge field in physics is denoted by $A_i$, math notation is $\nabla: E \to E \otimes \Omega^1A$ and sense of math notation is strictly definite. I think that written by physicists texts are accessible to understand for people which have contacts with physicists. I had no contacts with physicists, and so I did not understand their texts. Moreover notation of noncommutative geometry is much more clear then notation of commutative one. For example Levi-Civita connection is given by Christoffel symbols
\begin{equation*}
\Gamma^l_{jk}=\frac{1}{2}\sum_{r}g^{lr}\left\{\partial_k g_{ri} + \partial_{j}g_{rk}-\partial_r g_{jk}\right\}.
\end{equation*}
But noncommutative notation \cite{varilly:noncom}
\begin{equation*}
\left(r, \nabla s\right) - \left(\nabla r,s\right) = [D, (r|s)]
\end{equation*}
is much more clear.
Written by physicists texts contains a lot of heuristic formulas and it is difficult to distinguish them from rigorous ones. I did not find articles devoted to noncommutative Wilson lines which use math notation. So I have used both physical notation and math one. The correspondence between both these notations is explained.  This article has rather math style, so I assigned the mathematics category to it. 

\section{Introduction}
\paragraph{}In the commutative gauge theory the curvature of a connection is not sufficient to extract the complete gauge invariant information
 of the connection. This is the well-known gauge copy problem.
For certain gauge groups, the gauge copy problem can be
  solved by considering the set of all Wilson loops as
  the basic observables  of the gauge theory, see e.g.~\cite{Driver,Sengupta}.
 There are examples showing that there is also a gauge copy problem in noncommutative geometry,
 see e.g.~Proposition 4.2 from \cite{schenkel:nc_parallel}.
 This means that considering only observables which are constructed from the curvature, we can in general not extract the complete
 gauge invariant information of the connection. This calls for a suitable generalization of parallel transport and Wilson loops
 to noncommutative geometry. Some examples of noncommutative Wilson lines are described in \cite{alekseev_bytsko:wilson_nc_tori,Ambjorn:2000cs,Ishibashi:1999hs}. However these articles contain rather particular examples than a general theory. These examples can be used for the noncommutative torus only. General theory requires a noncommutative generalization of closed paths. However closed paths can be replaced with covering projections. A generalization of covering projections is described in my articles \cite{ivankov:cov_pr_nc_torus_fin,ivankov:infinite_cov_pr}. Following table contains necessary ingredients and their noncommutative analogues.
 \newline
 \break
 \begin{tabular}{|c|c|}
 \hline
 Differential geometry & Noncommutative geometry\\
 \hline
Spin manifold \cite{hajac:toknotes} & Spectral triple \cite{varilly:noncom} \\
 Vector bundle \cite{kobayashi_nomizu:diff_geom} & Projective module \cite{varilly:noncom} \\
 Gauge field \cite{green_schwarz_witten:superstring} & Noncommutative connection \cite{connes:c_alg_dg,hajac:toknotes} \\
 Parallel transport \cite{kobayashi_nomizu:diff_geom}  & Noncommutative parallel transport \cite{connes:c_alg_dg,schenkel:nc_parallel} \\
 Closed path \cite{spanier:at}\  & Noncommutative covering projection \cite{ivankov:cov_pr_nc_torus_fin,ivankov:infinite_cov_pr} \\

 \hline
 \end{tabular}
 \newline
 \newline
 \break
 Composition of these ingredients supplies a noncommutative generalization of Wilson lines.

\paragraph{} Following notation is used in this article.

 \begin{tabular}{|c|c|}
 \hline
 Symbol & Meaning\\
 \hline
 $\mathrm{Aut}(A)$ & Group * - automorphisms of $C^*$-algebra $A$\\
 $\mathbb{C}$ (resp. $\mathbb{R}$)  & Field of complex (resp. real) numbers \\
$C(X)$ & $C^*$ - algebra of continuous complex valued \\
 & functions on compact topological space $X$\\
$C_b(X)$ & $C^*$ - algebra of bounded  continuous complex valued \\
 & functions on locally compact topological space $X$\\
$C_0(X)$ & $C^*$ - algebra of continuous complex valued \\
 & functions on locally compact topological space which tends to 0 at infinity $X$\\

$M(A)$  & Multiplier algebra of $C^*-$ algebra $A$\\
 $\mathbb{N}$ & The set of natural numbers\\
$\mathbb{Q}$ & The field of rational numbers\\

 $U(A) \in A $ & Group of unitary operators of algebra $A$\\

 $\mathbb{Z}$ & Ring of integers \\

 $\mathbb{Z}_m$ & Ring of integers modulo $m$ \\

 \hline
 \end{tabular}
 \section{Prototype. Commutative Wilson Lines}\label{comm_prot}
 \begin{empt}{\it Standard description of Wilson lines} \cite{green_schwarz_witten:superstring}. Here the physical notation is used.
 Let $M$ be a compact Riemann manifold, and let $A$ be a locally pure gauge field, i.e. $A$ can be locally represented by following way
 \begin{equation}\label{pure_gauge}
 A_i = \partial_i U \cdot U^{-1}
 \end{equation}
 where $U$ is the gauge transformation. Whether \eqref{pure_gauge} is also true globally is another story. Condition \eqref{pure_gauge} is equivalent to that the field stength equals to zero \cite{green_schwarz_witten:superstring}.
 If the manifold is not simply connected, $\pi_1(M) \neq \{e\}$ then there is a more general possibility that appears in electrodynamics as the Bohm-Aharonov effect. Let $\gamma$ be a noncontractible loop in $M$, beginning and ending at the same point $x$. Then "Wilson line"
 \begin{equation}\label{gauge_transformation}
 U_{\gamma} = P \ \mathrm{exp} \oint\limits_{\gamma} A \cdot dx
 \end{equation}
 is gauge invariant, and if $U_{\gamma}\neq 1$ it cannot be set to one by gauge transformation, $U_{\gamma}$ depends only on $[\gamma] \in \pi_1(M)$. If $G$ is the gauge transformation group then is a group homomorphism
 \begin{equation}\label{from_fund_to_gauge}
 \varphi: \pi_1(M)\to G,
 \end{equation}
 \begin{equation*}
 [\gamma] \mapsto U_{\gamma}.
 \end{equation*}
 \end{empt}
 \begin{empt}\label{alter_prot} {\it Alternative description of Wilson lines}.
 In this construction I follow to \cite{green_schwarz_witten:superstring} 16.4.1. Let $\widetilde{M}$ be a simply connected manifold, and let $F$ be a discrete symmetry group which acts freely. Let $M = \widetilde{M}/F$. An ordinary field $\psi(x)$ is equivalent to a field on $\widetilde{M}$ that obeys
 \begin{equation}\label{invariant_simple}
 \psi(fx)=\psi(x); \ \forall f \in F.
 \end{equation}
 Then we can generalize \eqref{invariant_simple} as follows. There is a natural isomorphism $\pi_1(M) \approx F$, and from \eqref{from_fund_to_gauge} it follows that there is a natural homomorphism
 \begin{equation}\label{from_trans_to_field}
 F\to G,
 \end{equation}
 \begin{equation*}
 f \mapsto U_f.
 \end{equation*}
 Now we require that $\psi$ obey not \eqref{invariant_simple} but
 \begin{equation}\label{invariant_twist}
 \psi(fx)=U_f\psi(x); \ \forall f \in F.
 \end{equation}
 This operation enables us replace guage field \eqref{gauge_transformation} by "twist" in boundary conditions \eqref{invariant_twist}. The gauge field $A$ which obeys \eqref{gauge_transformation} (resp. field $\psi$) is replaced with the trivial gauge field (resp. field $\psi'$ which obey \eqref{gauge_transformation}). Field $\psi'$ is given by
 \begin{equation*}
 \psi'(y) = \left(P \ \mathrm{exp} \int\limits_{\omega} A' \cdot dx\right) \psi(\pi(y))
 \end{equation*}
 where
 \begin{itemize}
 \item $\pi: \widetilde{M} \to M$ is a covering projection.
 \item $\omega : [0,1] \to \widetilde{M}$ is such that $\omega(0) = y_0$ is a fixed point, $\omega(1)=y$.
 \item $A'$ is a lift of $A$ by $\pi$.
 \end{itemize}
 This construction is similar to passive/active approach to physical transformations. Change of a gauge field is similar to a passive transformation of a coordinat system, a change of field is similar to an active transformation of point's position. Both transformations describe the same physical phenomenon.
 \end{empt}

 \begin{empt}{\it Replacement of closed paths by covering projections}.
  Construction from \cite{green_schwarz_witten:superstring} can be generalized such that the universal covering is replaced by a covering  $\pi:\widetilde{M} \to M$ such that the image of the composition
  \begin{equation*}
  \pi_1(\widetilde{M})\xrightarrow{\pi_1(\pi)} \pi_1(M) \xrightarrow{\varphi} G
  \end{equation*}
  coincides with image of $\varphi$. In this case homomorphism \eqref{from_fund_to_gauge} can be replaced with
  \begin{equation}\label{from_trans_to_gauge}
  G(\widetilde{M}|M) \to G
  \end{equation}
  where $G(\widetilde{M}|M)$ is a group of covering transformations \cite{spanier:at}.
 \end{empt}
\section{Noncommutative Parallel Transport}\label{sec:fuzzy}
\subsection{Parallel transports}
\begin{empt} Definitions of section \ref{comm_prot} cannot be directly used in the noncommutative case because the noncommutative geometry does contain closed paths. However paths can be replaced with module parallel transports.
\end{empt}
\begin{defn}\label{defi:transport}\cite{schenkel:nc_parallel}
Let $A$ be an associative and unital algebra and $\mathcal{E}$ a right $A$-module.
\begin{itemize}
\item[1.)]
A {\it one-parameter group of automorphisms} of $A$ is a map
$\varphi: \mathbb{R} \times A \to A\,,~(\tau,a) \mapsto \varphi(\tau,a) =\varphi_\tau(a)$, such that
\begin{itemize}
\item[(i)] $\varphi_\tau(a\,b) = \varphi_\tau(a) \,\varphi_\tau(b)$, for all $\tau\in\mathbb{R}$ and $a,b\in A$
\item[(ii)] $\varphi_0 =\mathrm{id}_A$
\item[(iii)] $\varphi_{\tau+\sigma} = \varphi_\tau\circ\varphi_\sigma$, for all $\tau,\sigma\in\mathbb{R}$
\end{itemize}
\item[2.)] Let $\varphi:\mathbb{R} \times A \to A $ be a one-parameter group of automorphisms of $A$. A {\it module parallel transport} on
$\mathcal{E}$ along $\varphi$ is a map $\Phi: \mathbb{R} \times \mathcal{E} \to \mathcal{E}\,,~(\tau,s)\mapsto \Phi(\tau,s)=\Phi_\tau(s)$, such that
\begin{itemize}
\item[(i)] $\Phi_\tau(s\,a) = \Phi_\tau(s)\,\varphi_\tau(a)$, for all $\tau\in\mathbb{R}$, $s\in\mathcal{E}$ and $a\in A$
\item[(ii)] $\Phi_0=\mathrm{id}_\mathcal{E}$
\item[(iii)] $\Phi_{\tau+\sigma}= \Phi_\tau\circ \Phi_\sigma$, for all $\tau,\sigma\in\mathbb{R}$.
\end{itemize}
\end{itemize}
\paragraph{} If $A$ and $\mathcal{E}$ are equipped with a smooth structure, the maps $\varphi$ and $\Phi$ are required to be smooth. Denote by $\mathrm{Trans}_{\mathcal{E}}$ the set of module parallel transports.
\end{defn}
\begin{empt}
There are different notions of connections \cite{connes:c_alg_dg,parta:two_approaches_ym,schenkel:nc_parallel}, some of them are compared in \cite{parta:two_approaches_ym}.  The space of all connections on $\mathcal{E}$ is denoted by $\mathrm{Con}_A(\mathcal{E})$.
\end{empt}
\begin{defn}
Let $A$ be a $\mathbb{C}$-algebra, and let $\mathcal{E}$ be a finite projective $A$ module. Suppose that $\mathrm{Paths}_A$ is a set of one-parameter group of $A$ automorphisms.   A {\it connection transport procedure} is a natural map
\begin{equation*}
\mathbf{Transport} : \mathrm{Con}_A(\mathcal{E}) \times \mathrm{Paths}_A  \to \mathrm{Trans}_{\mathcal{E}}
\end{equation*}
such that
$\mathbf{Transport}(\nabla, \varphi)$ is a transport along $\varphi$  for any $\nabla \in \mathrm{Con}_A(\mathcal{E})$ and $\varphi \in \mathrm{Paths}_A$.
\end{defn}
\begin{rem}
A connection transport procedure should have a good math and/or physical sense. Such procedures are known in following cases:
\begin{enumerate}
\item Commutative differential geometry \cite{kobayashi_nomizu:diff_geom,schenkel:nc_parallel}.
\item $A = \mathbb{M}_n(\mathbb{C})$ \cite{schenkel:nc_parallel}.
\item Noncommutative torus \cite{alekseev_bytsko:wilson_nc_tori,Ambjorn:2000cs,Driver,Sengupta}.
\item Noncommutative differential geometry \cite{connes:c_alg_dg}.
\end{enumerate}
\end{rem}
\subsection{Parallel transform in noncommutative differential geometry}
\paragraph{} In this section I follow to \cite{connes:c_alg_dg}.
\begin{empt}\cite{connes:c_alg_dg}
Let $(A,G,\alpha)$ be a $C^*$ dynamical system, where $G$ is a Lie group.
We shall say that $x \in A$ is of $C^{\infty}$ class iff
the map $g \mapsto \alpha_g (x)$ from $G$ to the normed space $A$ is in
 $C^{\infty}$. The involutive algebra $A^{\infty} = \{ x \in A ,
\, x \ \hbox{of class} \ C^{\infty} \}$ is norm dense in $A$.

Let $\mathcal{E}^{\infty}$ be a finite projective module on $A^{\infty}$,
 (we shall write it as a right module); $\mathcal{E} = \mathcal{E}^{\infty}
\otimes_{A^{\infty}}$ $A$ is then  a finite projective module on
$A$.
\end{empt}

\begin{lem}\cite{connes:c_alg_dg}
For every finite projective module $\mathcal{E}$ on $A$, there exists a
finite projective module $\mathcal{E}^{\infty}$ on $A^{\infty}$, unique up to
isomorphism, such that $\mathcal{E}$ is isomorphic to $\mathcal{E}^{\infty}
\otimes_{A^{\infty}} \, A$.
\end{lem}
\begin{empt}
Let $\delta$ be the representation of Lie$G$ in the Lie-algebra of derivations of
 $A^{\infty}$ given by
$$
\delta_X (x) = \lim_{t \rightarrow 0} \ \frac{1}{t} \, (\alpha_{g_t} (x) - x) \, ,
\qquad \hbox{where} \quad \dot{g}_0 = X \, , \ x \in A^{\infty} \, .
$$
\end{empt}
\medskip

\begin{defn}\label{connes_dg_connection_defn}\cite{connes:c_alg_dg} $\mathcal{E}^{\infty}$ be a finite projective module on $A^{\infty}$,
a {\it connection} (on $\mathcal{E}^{\infty}$) is a linear map
$\nabla : \mathcal{E}^{\infty} \to \mathcal{E}^{\infty} \otimes(\mathrm{Lie} \ G)^*$ such that,
for all $X \in \mathrm{Lie} \ G$ and $\xi \in \mathcal{E}^{\infty}$, $x \in
A^{\infty}$ one has
$$
\nabla_X (\xi \cdot x) = \nabla_X (\xi) \cdot x + \xi \cdot \delta_X (x) \, .
$$
\end{defn}

\begin{empt}\label{connes_transp_procedure}
There is a natural correspondence between elements of Lie algebra and one-parameter transformation groups \cite{kobayashi_nomizu:diff_geom}.
Let $\nabla$ be a connection (on $\mathcal{E}^{\infty}$). If $X\in \mathrm{Lie}G$  defines a one-parameter group of automorphisms $\varphi$ of $A^{\infty}$  then $\nabla_X$ defines a  module parallel transport $\Phi$ on
$\mathcal{E}^{\infty}$ along $\varphi$. So we have a connection transport procedure
\begin{equation*}
\mathbf{Transport} : \mathrm{Con}_{A^{\infty}}(\mathcal{E}^{\infty}) \times \mathrm{Paths}_A^{\infty}  \to \mathrm{Trans}_{\mathcal{E}^{\infty}}.
\end{equation*}

\end{empt}

\begin{defn}\cite{connes:c_alg_dg}
Let $\nabla$ be a connection on the finite projective module
$\mathcal{E}^{\infty}$ (on $A^{\infty}$), the curvature of $\nabla$
is the element
$\mathcal{T}$ of $\mathrm{End}_{A^{\infty}} (\mathcal{E}^{\infty}) \otimes\mathbb{L}^2 (\mathrm{Lie} \
G)^*$ given by
$$
\mathcal{T} (X,Y) = \nabla_X \nabla_Y - \nabla_Y \nabla_X - \nabla_{[X,Y]} \in \mathrm{End}_{A^{\infty}} (\mathcal{E}^{\infty}) \, , \qquad \forall \, X,Y \in \mathrm{Lie} \ G
\, .
$$
\end{defn}
\begin{defn}\cite{connes:c_alg_dg} A connection with zero curvature is said to be {\it flat}.
Denote by $\mathrm{Con}_A(\mathcal{E})_0$ a space of flat connections.
\end{defn}

\begin{empt}\label{connection_construction}\cite{connes:c_alg_dg}
If $e \in A^{\infty}$ is an idempotent then every connection on $e \, A^{\infty}$ is of the form $\nabla_X (\xi) =
\nabla_X^0 (\xi) + \theta_X \, \xi$, $\forall \, \xi \in e \, A^{\infty}$, $X \in
\hbox{Lie} \ G$, where the form $\theta \in e \, \Om^1 \, e$ is uniquely
determined by $\nabla$, one has $\theta_X^* = -\theta_X$, $\forall \, X \in
\hbox{Lie}
\ G$ iff $\nabla$ is compatible with the hermitian structure of $e \,
A^{\infty}$. We identify $\mathrm{End} (e \, A^{\infty})$ with $e \, A^{\infty} \, e \subset
A^{\infty}$, the curvature $\mathcal{T}_0$ of the grassmannian connection is
the 2-form $e(de \wedge de) \in \Om^2$, the curvature of $\nabla = \nabla^0 +
\theta \wedge$ equals to
\begin{equation}\label{curvature_formula}
\mathcal{T}_0 + e (d\theta + \theta \wedge \theta) \, e \in \Om^2.
\end{equation}

\end{empt}

\section{Noncommutative Generalization of Closed Paths}
\begin{empt}{\it Commutative loops}. Some of parallel transports can be regarded as noncommutative loops. First of all we consider the commutative case. Let $M$ be a manifold and let $\pi: \widetilde{M} \to M$ be a covering projection, $G(\widetilde{M}|M)$ is the group of covering transformations, i.e. $M \approx \widetilde{M}/G(\widetilde{M}|M)$. Let $\varphi : \mathbb{R} \times C^{\infty}(M)\to C^{\infty}(M)$ satisfies condition 1 of definition \ref{defi:transport}. From Gelfand - Na\u{i}mark theorem \cite{arveson:c_alg_invt} it follows that $\varphi$ defines a one-parameter group of homeomorphisms (indeed diffeomorphisms) $\varphi^*: \mathbb{R} \otimes M \to M$. If $\varphi_1 = \mathrm{Id}_{C^{\infty}(M)}$ then a function $\gamma: [0, 1] \to M$, $t \mapsto \varphi^*(t, x_0)$ corresponds to a closed path. We would like to know whether this path is not contractible. It is known \cite{spanier:at} that covering projections have property of unique path lifting. So $\varphi^*$ can be lifted to $\widetilde{\varphi}^*: \mathbb{R} \times C^{\infty}(\widetilde{M})\to C^{\infty}(\widetilde{M})$. If $\widetilde{\varphi}^*_1 \in G(\widetilde{M}|M)$ then $\varphi^*_1 = \mathrm{Id}_M$. If $\widetilde{\varphi}^*_1$ is not a trivial element in $G(\widetilde{M}|M)$ then $\varphi$ does not correspond to a contractible path.
\end{empt}
\begin{empt}{\it Noncommutative loops.}
In my articles \cite{ivankov:cov_pr_nc_torus_fin,ivankov:infinite_cov_pr} there is a construction of covering projections for $C^*$-algebras (See Appendix \ref{appendix_nc_cov}). A  covering projection of $C^*$-algebra $A$ is a *-homomorphism $A \to M(\widetilde{A})$ where $\widetilde{A}$ is another $C^*$-algebra. A group of noncommutative covering transformations $G(\widetilde{A}|A)$ acts on $\widetilde{A}$ such that $g(a\widetilde{a})=ag\widetilde{a}$; $\forall a\in A$, $\forall \widetilde{a} \in \widetilde{A}$ $\forall g \in G(\widetilde{A}|A)$. For any $\varphi: \mathbb{R} \times A \to A$ there is an unique lift $\widetilde{\varphi}: \mathbb{R} \times \widetilde{A}\to\widetilde{A}$ such that
\begin{equation*}
\widetilde{\varphi}_{\tau}(a \widetilde{a})=\varphi_{\tau}(a)\widetilde{\varphi}_{\tau}(\widetilde{a}); \ a\in A, \ \widetilde{a} \in \widetilde{A}, \ \tau \in \mathbb{R}.
\end{equation*}
\begin{defn}
If $\varphi$ is such that
\begin{equation*}
\widetilde{\varphi}_1 = g \in G(\widetilde{A}|A),
\end{equation*}
\begin{equation*}
\widetilde{\varphi}_{\tau} \notin G(\widetilde{A}|A), \ 0 < \tau < 1.
\end{equation*}
then we say that $\varphi$ is a {\it closed path associated with $g$}.
\end{defn}

\end{empt}
\begin{defn}
Suppose that for any $g \in G(\widetilde{A}|A)$, $\nabla \in \mathrm{Con}_{A}(\mathcal{E})_0$ and for any closed paths $\varphi', \varphi''$ associated with $g$ we have
\begin{equation*}
\mathbf{Transport}(\nabla, \varphi'_1)=\mathbf{Transport}(\nabla, \varphi''_1).
\end{equation*}
A {\it generalized Wilson line} is a map
\begin{equation*}
\mathbf{Wilson} : G(\widetilde{\mathcal{A}}|\mathcal{A}) \times \mathrm{Con}_{\mathcal{A}}(\mathcal{E})_0 \to \mathrm{Aut}(\mathcal{E}),
\end{equation*}
\begin{equation*}
(g, \nabla) \mapsto \mathbf{Transport}(\nabla, \varphi_1)
\end{equation*}
where $\varphi$ is associated with $g$ and $\mathbf{Transport}$ is defined in \ref{connes_transp_procedure}.

\end{defn}

\section{Wilson Lines and Noncommutative Covering Spaces}
\paragraph{} Commutative geometry has a lot of local structures, for example local sections of bundles. There are  bundles such that they have no global sections. Let $p: P \to M$ is a bundle such that $p$ has no global sections. Let $\pi : \widetilde{M} \to M$ be the universal covering projection, and $\widetilde{p}: \widetilde{M} \to M$ be the pullback \cite{spanier:at} of $p$ by $\pi$. Then $\widetilde{p}$ can have global sections. Any locally pure gauge field can be locally represented by \eqref{pure_gauge}, but cannot be represented by \eqref{pure_gauge} globally in general case. However pullback of this field can be globally repesented  by \eqref{pure_gauge}. Noncommutative geometry has no local sections. However there is the noncommutative generalisation of covering projections \cite{ivankov:cov_pr_nc_torus_fin,ivankov:infinite_cov_pr}. Local pure gauge fields can be regarded as global pure gauge fields on universal noncommutative covering projections. Locally gauge fields satisfy \eqref{pure_gauge} which can be rewritten in physical notation \cite{alekseev_bytsko:wilson_nc_tori}
\begin{equation} \label{Ag}
  \partial_i g = i A_i * g \,.
\end{equation}
This equation has a noncommutative analog.  Let $A^{\infty}$ be a smooth algebra and $\mathcal{E}^{\infty}$ be a finitely generated projective smooth $A^{\infty}$ module.
Let $(A^{\infty}, \widetilde{A^{\infty}}, G, _{A^{\infty}}M_{\widetilde{A^{\infty}}})$ be a noncommutative covering projection (See Appendix \ref{appendix_nc_cov}). Spaces $A^{\infty}$,  $\widetilde{A^{\infty}}$, $\mathcal{E}^{\infty}$ are operator spaces. Let denote $\widetilde{\mathcal{E}^{\infty}}=\mathcal{E}^{\infty} \otimes _{A^{\infty}}\widetilde{A^{\infty}}$ where $\otimes$ means the Haagerup tensor product. The $\widetilde{\mathcal{E}^{\infty}}$ module can be regarded as pullback of $\mathcal{E}^{\infty}$. Any vector $X \in \mathrm{Lie}G$ can be lifted to the vector $\widetilde{X}\in \mathrm{Lie}\widetilde{G}$ where $\mathrm{Lie}\widetilde{G}$ is a Lie algebra of infinitesimal transformations of the $\widetilde{A^{\infty}}$. Any connection $\nabla : \mathcal{E}^{\infty} \to \mathcal{E}^{\infty} \otimes \Omega^1A^{\infty}$ can be lifted to $G$ equivariant connection $\widetilde{\nabla}: \widetilde{\mathcal{E}^{\infty}} \to \widetilde{\mathcal{E}^{\infty}} \otimes \Omega^1\widetilde{A^{\infty}}$. Otherwise any connection $\widetilde{\nabla}: \widetilde{\mathcal{E}^{\infty}} \to \widetilde{\mathcal{E}^{\infty}} \otimes \Omega^1\widetilde{A^{\infty}}$ naturally induces a map $\widetilde{\nabla}':\mathrm{End}_{\widetilde{A^{\infty}}}\left(\widetilde{\mathcal{E}^{\infty}}\right) \to \mathrm{End}_{\widetilde{A^{\infty}}}\left(\widetilde{\mathcal{E}^{\infty}}\right) \otimes \Omega^1\widetilde{A^{\infty}}$.
Noncummutative analog of \eqref{Ag} is given by
\begin{equation}\label{nc_Ag}
\widetilde{X}U = \widetilde{\nabla}_{\widetilde{X}}U; \ \forall X \in \mathrm{Lie}G, \ U\in \mathrm{Aut} \left(\widetilde{\mathcal{E}^{\infty}}\right).
\end{equation}
$U$ is a noncommutative analog of a global gauge transformation on the universal covering space.
Wilson line can be regarded as a group homomorphisms $G \to \mathrm{Aut}\left(\widetilde{\mathcal{E}^{\infty}}\right)$ given by
\begin{equation}\label{wilson_rel}
g \mapsto (gU) \cdot U^{-1}
\end{equation}
where $g \in G(\widetilde{A^{\infty}}|A^{\infty})$ and $U$ satisfies \eqref{nc_Ag}.

\section{Alternative Description of  Wilson Lines}
\paragraph{} As well as in \ref{alter_prot} (See \cite{green_schwarz_witten:superstring}) we can define an alternative approach to Wilson lines. Suppose that there is a spectral triple $(\mathcal{A}, H, D)$ \cite{varilly:noncom}. In \cite{ivankov:cov_pr_nc_torus_fin,ivankov:infinite_cov_pr} I defined noncomutative covering projection $(\widetilde{A}, \widetilde{H}, \widetilde{D})$ with a group of covering transformations $G(\widetilde{\mathcal{A}}|\mathcal{A})$, such that there is the natural *-homomorphism $\mathcal{A} \to M(\widetilde{\mathcal{A}})$ and representation $\pi: \mathcal{A} \to U(\widetilde{H})$.
Suppose that $\rho: G(\widetilde{\mathcal{A}}|\mathcal{A}) \to U(\widetilde{H})$ is a representation such that
\begin{equation*}
\widetilde{D} \rho(g)h = \rho(g) \widetilde{D}h, \ \forall h \in \mathrm{Dom}(\widetilde{D}).
\end{equation*}
Let $\widetilde{H}_0 = \{h \in \widetilde{H} \ | \ \pi(g)h = \rho(g)h\}$. Algebra $A$ naturally acts on $\widetilde{H}_0$ and we have a twisted spectral triple $(\mathcal{A}, \widetilde{H}_0, \widetilde{D})$ which can be regarded as description of Wilson line.

\section{Wilson Lines on the Noncommutative Torus}

\begin{empt}\label{sample_flat_conn}{\it A flat connection.} Let $\theta \in [0,1] - \mathbb{Q}$ and $A_{\theta}$ the $C^*$-algebra generated by two unitaries $u$, $v$ such that $u \, v = \lb \, v \, u$, $\lb = \exp ( 2 \,
\pi i \, \theta)$ which is said to be a noncommutative torus. There is a pre-$C^*$-algebra $\mathcal{A}_\theta$ of smooth functions defined in \cite{varilly:noncom}. There are  one-parameter groups $\varphi^{u}$, $\varphi^{v}$ such that
\begin{equation}\label{vp_u}
\varphi^{u}_{\tau}(u)=e^{2\pi i\tau}u, \ \varphi^{u}_{\tau}(v) = v,
\end{equation}
\begin{equation}\label{vp_v}
\varphi^{v}_{\tau}(u)=u, \ \varphi^{v}_{\tau}(v) = e^{2 \pi i\tau} v.
\end{equation}
Let $\mathcal{E} = \mathcal{A}_{\theta}$.
Let $\omega = i (c_u du + c_v dv) \in \Omega^1\mathcal{A}_{\theta}$ be such that $c_u, c_v \in \mathbb{R}$ then $\omega_X^* = -\omega_X$, $\forall \, X \in
\hbox{Lie}
\ G$. From \ref{connection_construction} it follows that there is a connection  connection $\nabla$ given by
\begin{equation*}
a \mapsto a \otimes \omega; \ a \in \mathcal{E} = \mathcal{A}.
\end{equation*}
From \eqref{curvature_formula} it follows that curvature of $\nabla$ is equal to
\begin{equation*}
1(d1 \wedge d1) + d\omega + \omega \wedge \omega.
\end{equation*}
All summands of above equation equal to zero, so curvature of $\nabla$ is zero, i.e. $\nabla$ is flat. One parameter groups $\varphi^u$, $\varphi^v$ correspond to following module parallel transports.
\begin{equation*}
\Phi^u_{\tau}(a) = e^{ic_u\tau}a,
\end{equation*}
\begin{equation*}
\Phi^v_{\tau}(a) = e^{ic_v\tau}a.
\end{equation*}
\end{empt}
\begin{exm}
In my article \cite{ivankov:cov_pr_nc_torus_fin} I have found all noncommutative covering projections of $\mathcal{A}_{\theta}$ Let $\theta' = \theta/4$ and  $\mathcal{A}_{\theta'}$ is generated by $x, y$ such that $xy = e^{2\pi i \theta'}yx$. There is a *-homomorphism $\pi: \mathcal{A}_{\theta} \to \mathcal{A}_{\theta'}$ such that
\begin{equation}\label{two_listed}
\pi(u)=x^2, \ \pi(v)=y^2.
\end{equation}
The group of transformation coverings $G( \mathcal{A}_{\theta'}|\mathcal{A}_{\theta})$ is isomorphic to $\mathbb{Z}_2 \times \mathbb{Z}_2$ with two generators $g_u, g_v \in \mathbb{Z}_2 \times \mathbb{Z}_2$ such that
\begin{equation}
g_ux = - x, \ g_uy = y, \ g_vx = x, \ g_vy = -y.
\end{equation}
One-paramemeter groups $\varphi^u$ and $\varphi^v$ can be lifted to $\widetilde{\varphi}^u$ and $\widetilde{\varphi}^v$ such that
\begin{equation*}
\widetilde{\varphi}^u_{\tau}(x)=e^{\frac{2\pi i\tau}{2}}x, \ \widetilde{\varphi}^u_{\tau}(y) = y,
\end{equation*}
\begin{equation*}
\widetilde{\varphi}^v_{\tau}(x)=x, \ \widetilde{\varphi}^v_{\tau}(y) = e^{\frac{2\pi i\tau}{2}}y.
\end{equation*}
From above equations it follows that $\widetilde{\varphi}^u_1 = g_u$,  $\widetilde{\varphi}^v_1 = g_v$. So $\varphi^u$ and $\varphi^{v}$ are closed paths associated with $g_u$ and $g_v$ respectively.  Let $\nabla$ be a flat connection defined in \ref{sample_flat_conn}. Generalized Wilson line is given by
\begin{equation*}
\mathbf{Wilson}(g_u, \nabla)(a) = e^{2\pi ic_u}a,
\end{equation*}
\begin{equation*}
\mathbf{Wilson}(g_v, \nabla)(a) = e^{2\pi ic_v}a.
\end{equation*}
where $a \in \mathcal{E}$.

\paragraph{} Let $\widetilde{\mathcal{A}}_{\theta}$ be an infinte covering projection of noncommutative torus (See \ref{nt_cp}). It is known \cite{ivankov:cov_pr_nc_torus_fin} that the covering transformation group equals to $\mathbb{Z}^2$.
Let us consider an unitary operator $U \in M\left(\widetilde{\mathcal{A}}_{\theta}\right)$ given by
\begin{equation*}
U = \pi_u(\phi^u) \pi_v(\phi^v)
\end{equation*}
where functions $\phi^u, \phi^v \in C_b(\mathbb{R})$ are given by
\begin{equation*}
\phi^u(x)= e^{ic_ux}, \phi^v(x)= e^{ic_vx}, \ \forall x \in \mathbb{R}.
\end{equation*}
It is easy to show that $U$ satisfies to \eqref{nc_Ag}. A group homomorphism \eqref{wilson_rel} is given by
\begin{equation*}
n_1 \mapsto e^{2\pi i c_u}1_{M(\widetilde{\mathcal{A}}_{\theta})},  \ n_2 \mapsto e^{2\pi i c_v}1_{M(\widetilde{\mathcal{A}}_{\theta})}.
\end{equation*}
where $n_1,n_2 \in \mathbb{Z}^2$ are generators of the covering transformation group.
\end{exm}
\begin{empt}
Let $\mathcal{E} = \mathcal{A}^4_{\theta}$ be a free module and let $e_1,..., e_4 \in  \mathcal{E}$ be its generators. Let $\nabla: \mathcal{E} \to \mathcal{E} \otimes \Omega^1(\mathcal{A}_{\theta})$ be a connection given by
\begin{equation*}
\nabla e_1 = c_u e_2 \otimes du, \ \nabla e_2 = -c_u e_1 \otimes du, \ \nabla e_3 = c_v e_4 \otimes dv, \ \nabla e_4 = -c_v e_3 \otimes dv.
\end{equation*}
where $c_u, c_v \in \mathbb{R}$. Let $X, Y \in \mathrm{Lie}G$ correspond to one-parameters groups given by \ref{vp_u} and \ref{vp_v} respectively. We have $[X,Y]=0$ because given by \ref{vp_u} and \ref{vp_v} groups commute. A direct calculation shows that $\nabla_X\nabla_Y = \nabla_Y \nabla_X$, so we have
\begin{equation*}
\mathcal{T} (X,Y) = \nabla_X \nabla_Y - \nabla_Y \nabla_X - \nabla_{[X,Y]} =0.
\end{equation*}

 Since $X$ and $Y$ a generators of $\mathrm{Lie}G$ we have $\mathcal{T}\equiv 0$, i.e a connection $\nabla$ is flat.
\begin{exm}
Let $\pi: \mathcal{A}_{\theta} \to \mathcal{A}_{\theta'}$ be a noncommutative covering projection given by \eqref{two_listed}. Then generalised Wilson line is given by
\begin{equation*}
\mathbf{Wilson}(g_u, \nabla)(e) =    \begin{pmatrix}
\mathrm{cos}(2\pi c_u) 1_{M(\widetilde{\mathcal{A}}_{\theta})} & -\mathrm{sin}(2\pi c_u) 1_{M(\widetilde{\mathcal{A}}_{\theta})}  & 0 & 0\\ \mathrm{sin}(2\pi c_u) 1_{M(\widetilde{\mathcal{A}}_{\theta})} & \mathrm{cos}(2\pi c_u) 1_{M(\widetilde{\mathcal{A}}_{\theta})} & 0 & 0 \\
0 & 0 & 1_{M(\widetilde{\mathcal{A}}_{\theta})} & 0 \\
0 & 0 & 0 & 1_{M(\widetilde{\mathcal{A}}_{\theta})} \\

\end{pmatrix}e,
\end{equation*}
\begin{equation*}
\mathbf{Wilson}(g_v, \nabla)(e) =  =  \begin{pmatrix}
 1_{M(\widetilde{\mathcal{A}}_{\theta})} & 0 & 0 & 0\\
 0 & 1_{M(\widetilde{\mathcal{A}}_{\theta})} & 0 & 0\\
0 & 0 & \mathrm{cos}(2\pi c_v) 1_{M(\widetilde{\mathcal{A}}_{\theta})} & -\mathrm{sin}(2\pi c_v) 1_{M(\widetilde{\mathcal{A}}_{\theta})} \\ 0 & 0 & \mathrm{sin}(2\pi c_v) 1_{M(\widetilde{\mathcal{A}}_{\theta})} & \mathrm{cos}(2\pi c_v) 1_{M(\widetilde{\mathcal{A}}_{\theta})} \\
\end{pmatrix}e
\end{equation*}
where $e \in \mathcal{E}$.
Let $\mathbf{c}_u, \mathbf{s}_u, \mathbf{c}_v, \mathbf{s}_v \in C_b(\mathbb{R})$ be such that
\begin{equation*}
\mathbf{c}_u(x) = \mathrm{cos}(c_ux), \ \mathbf{s}_u(x) = \mathrm{sin}(c_ux), \ \mathbf{c}_v(x) = \mathrm{cos}(c_vx), \ \mathbf{s}_v(x) = \mathrm{sin}(c_vx); \ \forall x \in \mathbb{R}.
\end{equation*}
\paragraph{} Let $U \in \mathrm{Aut}_{\widetilde{A}_{\theta}}\left(\widetilde{\mathcal{E}}\right)$ be an unitary element  given by
\begin{equation*}
U  =   \begin{pmatrix}
\pi_u(\mathbf{c}_u) & -\pi_u(\mathbf{s}_u) & 0 & 0\\
\pi_u(\mathbf{s}_u) & \pi_u(\mathbf{c}_u) & 0 & 0\\
0 & 0 & \pi_v(\mathbf{c}_v) & -\pi_v(\mathbf{s}_v)\\
0 & 0 & \pi_v(\mathbf{s}_v) & \pi_v(\mathbf{c}_v)\\
\end{pmatrix}
\end{equation*}
where $\pi_u$, $\pi_v$ are defined in \ref{nt_cp}.
\newline
Element $U$ satisfies to \eqref{nc_Ag}. A group homomorphism \eqref{wilson_rel} is given by
\begin{equation*}
n_1 \mapsto \begin{pmatrix}
\mathrm{cos}(2\pi c_u) 1_{M(\widetilde{\mathcal{A}}_{\theta})} & -\mathrm{sin}(2\pi c_u) 1_{M(\widetilde{\mathcal{A}}_{\theta})}  & 0 & 0\\ \mathrm{sin}(2\pi c_u) 1_{M(\widetilde{\mathcal{A}}_{\theta})} & \mathrm{cos}(2\pi c_u) 1_{M(\widetilde{\mathcal{A}}_{\theta})} & 0 & 0 \\
0 & 0 & 1_{M(\widetilde{\mathcal{A}}_{\theta})} & 0 \\
0 & 0 & 0 & 1_{M(\widetilde{\mathcal{A}}_{\theta})} \\

\end{pmatrix},
\end{equation*}
\begin{equation*}
n_2 \mapsto \begin{pmatrix}
 1_{M(\widetilde{\mathcal{A}}_{\theta})} & 0 & 0 & 0\\
 0 & 1_{M(\widetilde{\mathcal{A}}_{\theta})} & 0 & 0\\
0 & 0 & \mathrm{cos}(2\pi c_v) 1_{M(\widetilde{\mathcal{A}}_{\theta})} & -\mathrm{sin}(2\pi c_v) 1_{M(\widetilde{\mathcal{A}}_{\theta})} \\ 0 & 0 & \mathrm{sin}(2\pi c_v) 1_{M(\widetilde{\mathcal{A}}_{\theta})} & \mathrm{cos}(2\pi c_v) 1_{M(\widetilde{\mathcal{A}}_{\theta})} \\
\end{pmatrix}.
\end{equation*}

\end{exm}
\end{empt}
\section{Appendix. Noncommutative covering projections}
\subsection{General Theory}\label{appendix_nc_cov}
Let us remind notion of noncommutative covering projection.
\begin{defn}\cite{ivankov:cov_pr_nc_torus_fin}

Let $_AX_B$ be a Hermitian $B$-rigged $A$-module, $G$ is finite or countable group such that
\begin{itemize}
\item $G$ acts on $A$ and $X$,
\item Action of $G$ is equivariant, i.e $g (a\xi) = (ga) (g\xi)$ , and $B$ invariant, i.e $g(\xi b)=(g\xi)b$ for any $\xi \in X$, $b \in B$, $a\in A$, $g \in G$,
\item Inner-product  of $G$ is equivariant, i.e $\langle g\xi, g \zeta\rangle_X = \langle\xi, \zeta\rangle_X$ for any $\xi, \zeta \in X$,  $g \in G$.
\end{itemize}
Then we say that  $_AX_B$ is a {\it $G$-equivariant $B$-rigged $A$-module}.
\end{defn}
If $B$, $A$, $G$, $_AX_B$ satisfy definition 5.6 \cite{ivankov:infinite_cov_pr} than we say that quadruple $(B, A, G, _AX_B)$ is an {\it infinite noncommutative covering projection}. Finite covering projections are particular cases of infinite ones. If $(B, A, G, _AX_B)$ is an infinite covering projection than $G$ acts on $A$ there is a *-homomorphism $\pi: B \to M(A)$ such that
\begin{equation}
g(\pi(b)a)=\pi(b)(ga), \ b \in B, \ a \in A, \ g \in G.
\end{equation}
If $\mathcal{B} \in B$ is a smooth subalgebra then there is a smooth version $(\mathcal{B}, \mathcal{A}, G, _{\mathcal{A}}\mathcal{X}_{\mathcal{B}})$ of covering projection defined in section 8 of \cite{ivankov:infinite_cov_pr}.
\subsection{Infinite Noncommutative Covering Projection of Noncommutative Torus}\label{nt_cp}
\paragraph{} Let $A_{\theta}$ be a noncommutative torus generated by unitary elements $u$, $v$ \cite{varilly:noncom}. In 7.2 \cite{ivankov:infinite_cov_pr} it is constructed an infinite covering projection $(A_{\theta}, \widetilde{A_{\theta}}, \mathbb{Z}^2, _{\widetilde{A}_{\theta}}X_{A_{\theta}})$. This construction contains a representation $A_{\theta} \rar B(H)$ and two representations $\pi_u : C_0(\mathbb{R}) \to B(H)$, $\pi_v : C_0(\mathbb{R}) \to B(H)$ such that $\widetilde{A}_{\theta}$ is the norm completion of generated by following operators
\begin{equation}\label{nc_ope}
\pi_u(f)\pi_v(g), \ \pi_v(f), \pi_u(g); \ f,g \in C_0(\mathbb{R}).
\end{equation}
subalgebra of $B(H)$.
\newline
Both $A_{\theta}$ and $\widetilde{A_{\theta}}$ are represented in the same Hilbert space $H$. These representations induce a *-homomorphism $A_{\theta}\to M(\widetilde{A}_{\theta})$ given by
\begin{equation*}
u \cdot \pi_u(f)\pi_v(g) = \pi_u(\varphi^{\mathrm{exp}}f)\pi_v(g),
\end{equation*}
\begin{equation*}
v \cdot \pi_u(f)\pi_v(g) = \pi_u(f)\pi_v(\varphi^{\mathrm{exp}}g).
\end{equation*}
where $\varphi^{\mathrm{exp}}\in C_b(\mathbb{R})$ be given by $x \mapsto e^{ix}$ ($x \in \mathbb{R}$).
For any $f\in C_0(\mathbb{R})$ denote by $f^{\uparrow}\in C_0(\mathbb{R})$ given by
\begin{equation*}
f^{\uparrow}(x)= f(x + 2\pi); \ \forall x \in \mathbb{R}.
\end{equation*}
Noncommutative group of covering transformation equals to $\mathbb{Z}^2$ and generators $n_1, n_2 \in \mathbb{Z}^2$ act on operators \eqref{nc_ope} by following way
\begin{equation*}
n_1 \cdot \pi_u(f)\pi_v(g) = \pi_u(f^{\uparrow}) \pi_v(g), \ n_2 \cdot \pi_u(f) \pi_v(g) = \pi_u(f) \pi_v(g^{\uparrow});
\end{equation*}
Extension of this action gives an action of $\mathbb{Z}^2$ on $\widetilde{A}_{\theta}$.

\end{document}